\documentclass[a4paper,abstract=true]{scrartcl}

\usepackage[utf8]{inputenc}
\usepackage{enumitem}

\usepackage{amssymb, amsmath}
\usepackage{bm}
\usepackage{mathtools}
\usepackage{ntheorem}

\usepackage{booktabs, array}

\usepackage{graphicx}
	\graphicspath{{graphics/}{experiments/denoising/}{experiments/denoising/noise/}{experiments/denoising/noise/3levels}{experiments/histograms/}{experiments/edges/}{experiments/estimation/}}

\usepackage{xcolor}
\definecolor{aaublue}{rgb}{0.00,0.41,0.66}
\definecolor{nicered}{rgb}{0.647,0.129,0.149}

\usepackage{hyperref}
\hypersetup{%
	breaklinks=true,%
	colorlinks=true,%
	linkcolor=aaublue,%
	citecolor=nicered,%
	urlcolor=aaublue,%
	pdftitle={Gaussian-log-Gaussian wavelet trees ...},%
	pdfauthor={Robert Dahl Jacobsen and Jesper M{\o}ller},%
	pdfsubject={},%
	pdfkeywords={},%
	pdfpagemode=UseOutlines,%
	pdfstartview=FitH%
}

\renewcommand{\v}[1]{\ensuremath{\bm{#1}}}

\newcommand{\m}[1]{\v{#1}}

\newcommand{\T}{\top}

\newcommand{\Beta}{\ensuremath{B}}
\newcommand{\Kappa}{\ensuremath{K}}

\DeclareMathOperator{\diag}{\text{diag}}

\newcommand{\R}{\ensuremath{\mathbb{R}}}

\newcommand{\E}{\ensuremath{\operatorname{E}}}
\newcommand{\cov}{\ensuremath{\operatorname{Cov}}}
\newcommand{\var}{\ensuremath{\operatorname{Var}}}

\renewcommand{\d}{\ensuremath{\mathrm{d}}}
\newcommand{\diff}{\ensuremath{\,\mathrm d}}

\renewcommand{\l}{r}

\DeclareMathOperator*{\argmax}{argmax}

\title{Frequentist and Bayesian inference for Gaussian-log-Gaussian wavelet trees and statistical signal processing applications}
\author{Robert Dahl Jacobsen\thanks{Department of Mathematical Sciences, Aalborg University, Skjernvej 4A, DK-9220 Aalborg East}
\and Jesper M{\o}ller\footnotemark[1]{}}

\begin{document}

\maketitle

\begin{abstract}
We introduce new estimation methods for a sub-class of the Gaussian scale mixture models for wavelet trees by Wainwright, Simoncelli \& Willsky that rely on modern results for composite likelihoods and approximate Bayesian inference.
Our methodology is illustrated for denoising and edge detection problems in two-dimensional images.
\end{abstract}


\section{Introduction}
\label{sec:intro}

Statistical models of wavelet coefficients in signal processing aims at modelling their non-Gaussian nature and statistical dependencies.
Different models have been proposed for this task, in particular a 
broad class of models for wavelet coefficients known as Gaussian scale mixtures 
\cite{crouse:nowak:baraniuk:98,Wainwright:Simoncelli:Willsky:2001}.
The present paper introduces new procedures for inference in a
subclass of these models that we refer to as
\emph{Gaussian-log-Gaussian wavelet tree models} or just GLG
models. 
Briefly, a GLG model specifies dependence in 
highpass wavelet coefficients with $d\ge1$ directional subbands
by a hidden Gaussian Markovian tree
structure such that coefficients
given the hidden states are independent and each coefficient 
follows a zero-mean Gaussian
distribution with the log-variance equal to the corresponding hidden
state. Further details of the GLG model are given in Section\nobreakspace \ref {sec:model}.

Section\nobreakspace \ref {sec:freq} introduces parameter estimation with moment matching
and an expectation-maximization (EM) algorithm. The moment-based estimates can be
used as initial values for the EM algorithm. The EM algorithm
is based on composite likelihoods and is fast and accurate if $d=1$
or if we assume independence for hidden states corresponding to different
directional subbands as is customary in many situations. 

Section\nobreakspace \ref {sec:bayes} demonstrates how the GLG model fits in the integrated
nested Laplace approximation (INLA)
framework \cite{rue:martino:07,rue:martino:chopin:09} for doing fast
approximative Bayesian inference. In particular the marginal posterior
distributions for the parameters can be estimated but the 
INLA approach requires a relatively low dimensional parameter as
compared with the dimenionality of the parameter in the full GLG
model. In a simulation study in Section\nobreakspace \ref {sec:simulations} we compare the results of the INLA approach with the results of the EM algorithm. 

As examples of applications Section\nobreakspace \ref {s:applications} discusses how the GLG model can be used for denoising and edge detection in images. 
 
Technical details are deferred to Appendices\nobreakspace  \ref{app:moments}-\ref{app:denoise} .
Matlab and R \cite{R} codes for our statistical inference procedures are available at \url{http://people.math.aau.dk/~robert/software}.

\section{The Gaussian-log-Gaussian wavelet tree model}
\label{sec:model}

This section recaps the GLG model \cite{Wainwright:Simoncelli:Willsky:2001} and sets the notation.

\begin{figure}
	\begin{center}
		\includegraphics{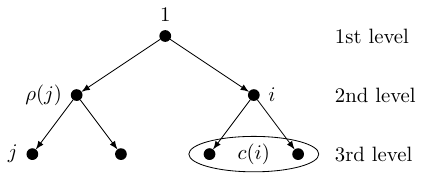}
	\end{center}
	\caption{Illustration of a binary tree structure corresponding
          to a one-dimensional signal (that is, $d=1$) 
with 3 levels of wavelet coefficients (that is, $L = 3$). 
Node $j$ has one parent $\rho(j)$ and node $i$ has two children $c(i)$.}
	\label{fig:binary_tree}
\end{figure}

The wavelet transform of a signal is separated into lowpass and
highpass coefficients. We wish to model the highpass coefficients, denoted
$\mathbf w = (\v w_1, \ldots, \v w_n)$, with respect to a given directed
tree, see
 Figure\nobreakspace \ref {fig:binary_tree} for a simple example. Here $i=1,\ldots,n$ 
is an abstract index 
identified with the nodes of the tree, the root is $i=1$ (the coarsest level
of the wavelet transform) and the directed 
edges correspond to the parent-child
relations of the 
 coefficients at the coarsest to the finest level. Specifically, we denote 
\begin{itemize}
\item $\l(i)$ the level of node $i$, that is, $\l(i)$ is the
number of nodes in the path from the root to $i$; thus $\l(1)=1$;
\item $L$ the number
of levels (the longest path);
\item
$c(i) \subset \{1,\ldots,n\}$ the children of $i$, that is, there is a
directed edge from $i$ to each child $j\in c(i)$; 
if $i$ is at the finest wavelet level, $i$ has no children ($c(i)=\emptyset$); 
else $c(i) \neq \emptyset$; 
in our image examples, $c(i)$ has cardinality 4 whenever
$c(i)\not=\emptyset$; and for convenience we let $c(0) = \{1\}$;
\item $\rho(j)$ the parent to node $j \neq 1$, that is, there is a
directed edge from $\rho(j)$ to $j$.
\end{itemize}
Furthermore, each 
$\v w_i = (w_i(1), \ldots, w_i(d))$ corresponds to
$d\ge 1$ directional subbands, where 
$d$ depends on the dimension of the signal. 
In our examples with images, $d = 3$.

We model the dependence structure for the coefficients through hidden
states as follows. 
Assume that each $\v w_i$ has an associated hidden state $\v s_i =
(s_i(1), \ldots, s_i(d))$ such that conditional on all hidden states
$\mathbf s = (\v s_1, \ldots, \v s_n)$,
the coefficients $w_i(\ell)$ ($\ell=1,\ldots,d$, $i=1,\ldots,n$) 
are independent and the conditional distribution of $\v w_i$
depends only on $\v s_i$ and has density
\begin{equation}
	p\bigl(\v w_i \,|\,\v s_i\bigr)
\sim N_d\bigl(\v 0, \diag\bigl[\exp(s_i(\ell)), \ell=1,\ldots,d\bigr]\bigr).
	\label{eq:obs_model}
\end{equation}
Here $N_d$ denotes the $d$-dimensional Gaussian distribution 
and the covariance matrix in 
\eqref{eq:obs_model} is diagonal.

We consider a parametric graphical model for the hidden states, 
using $\v\theta$ as generic
notation for the unknown parameters appearing in \textup {(\ref {eq:theta})}
later, and assume tying within levels as follows. 
Viewing $\mathbf s$ 
as a directed graphical model \cite{lauritzen:96}, we assume that the
conditional independence structure for $\v s_1, \ldots, \v s_n$ 
is given by the tree
structure. 
Hence the joint density
for coefficients and hidden states given $\v\theta$ has the structure
\begin{equation}
	\label{eq:full_model}
	p(\mathbf w, \mathbf s \,|\, \v\theta)
	= p\bigl(\v s_1 \,|\, \v\mu_{1}, \m\Sigma_{1}\bigr) \prod_{i=1}^n \Bigl[ p(\v w_i \,|\, \v s_i) \prod_{j\in c(i)} p\bigl(\v s_j \,|\, \v s_i, \v\alpha_{\l(i)}, \m\Beta_{\l(i)}, \m\Kappa_{\l(i)}\bigr) \Bigr],
\end{equation}
where we set the product over the empty set to be $1$ and where we
impose the
conditional densities
\begin{align}
	p\bigl(\v s_1 \,|\,\v\mu_{1}, \m\Sigma_{1}\bigr)
	& \sim N_d\bigl(\v\mu_{1}, \m\Sigma_{1}\bigr),
	\label{eq:root_dist}
	\\
	p\bigl(\v s_j \,|\,\v s_i, \v\alpha_{\l(i)}, \m\Beta_{\l(i)}, \m\Kappa_{\l(i)}\bigr)
	& \sim N_d\bigl(\v\alpha_{\l(i)} + \m\Beta_{\l(i)} \v s_i, \m\Kappa_{\l(i)}\bigr),
	\label{eq:transfer_dist}
\end{align}
Thus $\v\theta$
consists of 
the elements of the parameter vectors and matrices
\begin{gather}
	\begin{split}
		\label{eq:theta}
		&\v\mu_{1}  = [\mu_{1,\ell}]_{\ell=1,\ldots,d},
		\quad
		\m\Sigma_{1} = \bigl[\sigma_{1,\ell,\ell'}\bigr]_{\ell,\ell'=1,\ldots,d},
		\\
		&\v\alpha_{\l}  = [\alpha_{\l,\ell}]_{\ell=1,\ldots,d},
		\quad
		\m\Kappa_{\l} = \bigl[\kappa_{\l,\ell,\ell'}\bigr]_{\ell,\ell'=1,\ldots,d},
		\quad
		\m\Beta_{\l}  =
                \bigl[\beta_{\l,\ell,\ell'}\bigr]_{\ell,\ell'=1,\ldots,d},
	\end{split}
\end{gather}
with $\l=1,\ldots,L-1$.

We call \eqref{eq:full_model} the {\itshape (full) GLG 
(GLG) model} if we add no further restrictions than 
$\v\mu_{1}$ and $\v\alpha_{\l}$ being real $d$-dimensional
vectors, $\m B_{\l}$ being a  real $d\times d$ matrix, and 
 $\m\Sigma_{1}$ and $\m\Kappa_{\l}$ being symmetric and
strictly positive definite $d\times d$ matrices. The submodel
where the matrices 
$\m\Sigma_{1}$, $\m\Kappa_{1},\ldots,\m\Kappa_{L-1}$ and 
$\m\Beta_{1},\ldots,\m\Beta_{L-1}$ 
 are diagonal is called the {\itshape directional independence GLG model},
 because hidden variables corresponding to different directional subbands are
 independent. We also
 consider the simpler {\itshape homogenous GLG model} where for all
 $\ell=1,\ldots,d$ and $\l=1,\ldots,L-1$,  
\[\mu_{1,\ell}=\mu_1,\quad \m\Sigma_{1} =\sigma_1\m I_d,\quad 
\alpha_{\l,\ell}=\alpha,\quad \m\Kappa_{\l} =\kappa\m I_d,\quad
\m\Beta_{\l}=\beta\m I_d,\] 
where $\m I_d$ is the $d\times d$ identity matrix and $\mu_1\in\mathbb R$, 
$\sigma_1>0$, $\alpha\in\mathbb R$, $\kappa>0$ and $\beta\in\mathbb R$
are the five free parameters. The 
number of free parameters in the
full GLG model is $d + d(d+1)(3L/2 - 1)$,  that is, $18L - 9$ if $d=3$. 
 In the
directional independence GLG model we have $d(2+3L)$ free parameters,
that is, 
$9L+6$ if $d=3$. 

The inference procedures introduced later will be based on $k\ge1$
wavelet trees $\mathbf w^{(t)} = (\v w_1^{(t)}, \ldots, \v
w_n^{(t)})$, $t=1,\ldots,k$; 
in our examples, the pixels of an image correspond to the nodes of
the $k$ trees after applying
the (inverse) wavelet transform. We denote the corresponding hidden
states by $\mathbf s^{(t)} = (\v s_1^{(t)}, \ldots, \v
s_n^{(t)})$, $t=1,\ldots,k$, and assume that 
$(\mathbf w^{(1)},\mathbf s^{(1)}),\ldots,\mathbf w^{(k)},\mathbf
s^{(k)})$ are independent copies of $(\mathbf w,\mathbf s)$. We
denote all the coefficients (our data) and all the hidden states by 
\[{\bar{\mathbf w}}=(\mathbf w^{(1)},\ldots,\mathbf w^{(k)}),\quad
 {\bar{\mathbf s}}=(\mathbf s^{(1)},\ldots,\mathbf s^{(k)}).\]

\subsection{Mean and variance-covariance structure}

We shall later exploit the following marginal distributions for the
hidden states and the connection between their mean and
variance-covariance structure and the original parametrization in \eqref{eq:theta}.

By \textup {(\ref {eq:root_dist})} and\nobreakspace  \textup {(\ref {eq:transfer_dist})}, each hidden state $\v s_j$ follows a $d$-dimensional Gaussian distribution with a mean vector and covariance matrix depending only on the level $\l(j)$:
\begin{equation}
	\label{eq:marginal_node_pdf}
	p\left(\v s_j \ |\ \v\mu_{\l(j)}, \m\Sigma_{\l(j)}\right)
	\sim N_d\left(\v\mu_{\l(j)}, \m\Sigma_{\l(j)}\right)
\end{equation}
where the mean vector and covariance matrix
are determined recursively from the coarsest level to the second finest level by
\begin{equation}
	\label{eq:marginal_moment_relation}
	\v\mu_{r+1} = \v\alpha_r + \m\Beta_r \v\mu_r,
	\quad
	\m\Sigma_{r+1} = \m\Kappa_r + \m\Beta_r \m\Sigma_r
        \m\Beta_r^\T,
         \quad
         r=1,\ldots,L-1,
\end{equation}
where $\m\Beta_r^\T$ denotes the transpose of $\m\Beta_r$.
For a node $i$ at level $r=r(i)$ and
with children $j,h \in c(i)$, 
we obtain from \textup {(\ref {eq:transfer_dist})} that
\begin{align*}
	\cov(\v s_j, \v s_i) 
=  \m\Beta_{\l} \m\Sigma_{\l},\quad
	\cov(\v s_j, \v s_h) 
= \m\Beta_{\l} \m\Sigma_{\l} \m\Beta_{\l}^\T\quad\mbox{if $j\neq h$}.
\end{align*}


\section{Frequentist inference}
\label{sec:freq}

This section introduces estimating equations based on moment relations and an EM algorithm.
The main purpose of the moment method is to provide initial parameter estimates for the EM algorithm.

\subsection{Moment matching relations}
\label{sec:moment_relation}

Moments of the form $\E\bigl[w_i(\ell)^a w_j(\ell')^b\bigr]$ for $1\leq i,j \leq n$, $1\leq \ell, \ell'\leq d$ and $a,b \in \{0,1,\ldots\}$ can be derived by conditioning on the hidden states and exploiting well-known moment results for the log-Gaussian distribution.
Matching these expressions for the moments with empirical moments estimates allows us to estimate the parameters.
The details are deferred to Appendix\nobreakspace \ref {app:moments}.

\subsection{Composite likelihoods and the EM algorithm}
\label{sec:composite_EM}

The EM algorithm in Section~\ref{subsec:composite_EM} 
is based on composite likelihoods and to this end we need the marginal
distributions in Section~\ref{sec:marginal_likelihoods}.

\subsubsection{Marginal distributions for wavelet coefficients}
\label{sec:marginal_likelihoods}

First, combining \textup {(\ref {eq:obs_model})} and\nobreakspace  \textup {(\ref {eq:root_dist})}, we obtain the joint
density at the root, 
\begin{align}
	\MoveEqLeft
	p(\v s_1, \v w_1 \,|\, \v\mu_1, \m\Sigma_1) 
	= p(\v w_1 \,|\, \v s_1) p(\v s_1 \,|\, \v\mu_1, \m\Sigma_1) 
	\nonumber
	\\
	& = \frac1{(2\pi)^d \det(\m\Sigma_1)} \exp\biggl(-\frac12 \biggl[\biggl\{\sum_{\ell=1}^d \frac{w_1(\ell)^2}{\exp(s_1(\ell))} + s_1(\ell)\biggr\} + (\v s_1 - \v\mu_1)^\T \m\Sigma_1^{-1} (\v s_1 - \v\mu_1)\biggr]\biggr)
	\label{eq:jointpdf_root}
\end{align}
and thereby the marginal density of the root wavelet,
\begin{equation}
	\label{eq:marginal_root_pdf}
	p(\v w_1 \,|\, \v\mu_1, \m\Sigma_1) 
	= \int_{\R^d} p(\v s_1, \v w_1 \,|\, \v\mu_1, \m\Sigma_1) \d\v s_1.
\end{equation}
Hence the marginal log-likelihood based on the root wavelets $\overline{\v w}_1 = (\v w_1^{(1)}, \ldots, \v w_1^{(k)})$ for the $k$ trees is given by
\begin{equation}
	\label{eq:root_likelihood}
	l_1(\v\mu_1, \m\Sigma_1 | \overline{\v w}_1) 
	= \sum_{t=1}^k \log p(\v w_1^{(t)} \,|\, \v\mu_1, \m\Sigma_1).
\end{equation}

Second, consider any $i \in \{1,\ldots,n\}$ with level $r=\l(i) < L$.
Denote $\v w_{i,c(i)}$ the vector consisting of $\v w_i$ and all $\v w_j$ with $j\in c(i)$, and $\v s_{i,c(i)}$ the vector consisting of $\v s_i$ and all $\v s_j$ with $j\in c(i)$.
Using \textup {(\ref {eq:obs_model})},  \textup {(\ref {eq:transfer_dist})} and\nobreakspace  \textup {(\ref {eq:marginal_node_pdf})}, we obtain the density of $\bigl(\v s_{i,c(i)}, \v w_{i,c(i)}\bigr)$,
\begin{align}
	\MoveEqLeft[3]
	p(\v s_{i,c(i)}, \v w_{i,c(i)} \,|\, \v\mu_{r}, \m\Sigma_{r}, \v\alpha_r, \m\Beta_r, \m\Kappa_r) 
	\nonumber
	\\
	= {} & p(\v w_i \,|\, \v s_i) \ p(\v s_i \,|\, \v\mu_{r}, \m\Sigma_{r}) 
	\prod_{j\in c(i)} p(\v w_j \,|\, \v s_j) \ p(\v s_j \,|\, \v s_i, \v\alpha_r, \m\Beta_r, \m\Kappa_r)
	\nonumber
	\\
	= {} & \frac1{(2\pi)^{d(1+|c(i)|)} \det(\m\Sigma_{r})^{1/2} \det(\m\Kappa_r)^{|c(i)|/2}}
	\nonumber
	\\
	& \exp\biggl(-\frac12 \biggl[\biggl\{\sum_{\ell=1}^d \frac{w_i(\ell)^2}{\exp(s_i(\ell))} + s_i(\ell)\biggr\} + (\v s_i - \v\mu_{r})^\T \m\Sigma_{r}^{-1} (\v s_i - \v\mu_{r})
	\nonumber
	\\
	& + \sum_{j\in c(i)} \biggl\{\sum_{\ell=1}^d \frac{w_j(\ell)^2}{\exp(s_j(\ell))} + s_j(\ell)\biggr\} + (\v s_j - \v\alpha_r - \m\Beta_r \v s_i)^\T \m\Kappa_r^{-1} (\v s_j - \v\alpha_r - \m\Beta_r \v s_i)\biggr]\biggr),
	\label{eq:jointpdf_transfer}
\end{align}
where $|c(i)|$ denotes the number of children of node $i$, and
so  $\v w_{i,c(i)}$ has density
\begin{align}
	\MoveEqLeft
	p(\v w_{i,c(i)} \,|\, \v\mu_{r}, \m\Sigma_{r}, \v\alpha_r, \m\Beta_r, \m\Kappa_r) 
	\nonumber
	\\
	& = \int_{\R^{d(1+|c(i)|)}} 
	p(\v s_{i,c(i)}, \v w_{i,c(i)} \,|\, \v\mu_{r-1}, \m\Sigma_{r-1}, \v\alpha_r, \m\Beta_r, \m\Kappa_r) 
	\d\v s_{i,c(i)}.
	\label{eq:marginal_transfer_pdf}
\end{align}
Hence, denoting $\overline{\v w}_{i,c(i)}$ the vector of the $i$th wavelets $\v w_i^{(1)}, \ldots, \v w_i^{(k)}$ and their children $\v w_j^{(1)}, \ldots, \v w_j^{(k)}$, $j\in c(i)$, the log-likelihood based on $\overline{\v w}_{i,c(i)}$ is
\begin{equation}
	\label{eq:transfer_likelihood}
	l_r(\v\mu_{r}, \m\Sigma_{r}, \v\alpha_r, \m\Beta_r, \m\Kappa_r \,|\, \overline{\v w}_{i,c(i)}) 
	= \sum_{t=1}^k \sum_{j\in c(i)} p(\v w_r^{(t)}, \v w_j^{(t)} \,|\, \v\mu_{r}, \m\Sigma_{r}, \v\alpha_r, \m\Beta_r, \m\Kappa_r) .
\end{equation}

A common thread in the integrals
\textup {(\ref {eq:marginal_root_pdf})} and\nobreakspace  \textup {(\ref {eq:marginal_transfer_pdf})} is that
\eqref{eq:jointpdf_root} 
and \eqref{eq:jointpdf_transfer} are of the form
\begin{equation*}
	\text{non-linear function} \times \text{Gaussian density}.
\end{equation*}
This suggests that the Gauss-Hermite quadrature rule is a good choice
for evaluating the integral, see, for example, \cite{Press:Teukolsky:Vetterling:Flannery:2002}. Indeed,
under the directional independence GLG model or if $d = 1$, 
this is fast and accurate. However,
for the full GLG model with
$d = 3$, it is not feasible to calculate
\textup {(\ref {eq:marginal_transfer_pdf})} to a reasonable precision with
quadrature rules because of the large number of times it is required
in the EM algorithm in the following section---even with sparse grids as used in, for example,
\cite{Heiss:Winschel:2008}.

\subsubsection{The composite EM algorithm}
\label{subsec:composite_EM}

Because the following EM algorithm applies on composite likelihoods
\cite{Gao:Song:2011} defined from the marginal likelihoods in
Section\nobreakspace \ref {sec:marginal_likelihoods}, we call it the {\itshape composite EM algorithm}.
In brief, it proceeds from the coarsest to the finest level, using the
relation \textup {(\ref {eq:marginal_moment_relation})} for the parameters in the marginal and conditional
distributions of the hidden variables as follows.
\begin{enumerate}
	\item\label{item:root_EM}
		Apply the EM algorithm for the 
log-likelihood \textup {(\ref {eq:root_likelihood})} to obtain an estimate $(\widehat{\v\mu}_1, \widehat{\v\Sigma}_1)$.

	\item\label{item:transfer_EM} 
		For $r=1,\ldots,L-1$ proceed as follows. Denote
                $\overline{\v w}_{(r)}$ the vector of all
                $\overline{\v w}_{i,c(i)}$ with $\ell(i) = r$. Note
                that the log-composite likelihood given by the sum of the log-likelihoods \textup {(\ref {eq:transfer_likelihood})} based on all $\overline{\v w}_{i,c(i)}$ with $\ell(i) = r$ is 
	\begin{align*}\label{e:mmmm}
		\MoveEqLeft l_{r}(\v\mu_{r}, \m\Sigma_{r}, \v\alpha_r, \m\Beta_r, \m\Kappa_r \,|\, \overline{\v w}_{(r)}) 
		 = \sum_{i:\ell(i)=r}
		l_i(\v\mu_{r}, \m\Sigma_{r}, \v\alpha_r, \m\Beta_r, \m\Kappa_r \,|\, \overline{\v w}_{i,c(i)}).
	\end{align*}
	Then, using a previously obtained  estimate $(\widehat{\v\mu}_{r}, \widehat{\v\Sigma}_{r})$,
	apply the EM algorithm on  $l_{r}(\widehat{\v\mu}_{r}, \widehat{\m\Sigma}_{r}, \v\alpha_r, \m\Beta_r, \m\Kappa_r \,|\, \overline{\v w}_{(r)})$ to obtain an estimate $(\widehat{\v\alpha}_r, \widehat{\m\Beta}_r, \widehat{\v\Kappa}_r)$.
	Thereby, using \textup {(\ref {eq:marginal_moment_relation})}, an
        estimate $(\widehat{\v\mu}_{r+1}, \widehat{\m\Sigma}_{r+1})$ is also
        obtained. 
\end{enumerate} 
Appendix\nobreakspace \ref {app:compEM} details these steps.

\section{Bayesian inference}
\label{sec:bayes}

Section\nobreakspace \ref {sec:brief_INLA} provides a brief description of
approximate Bayesian inference for our GLG model using integrated
nested Laplace approximations.
Section\nobreakspace \ref {sec:simulations} compares the Bayesian approach with the EM algorithm of Section\nobreakspace \ref {subsec:composite_EM}.

\subsection{A short diversion into INLA}
\label{sec:brief_INLA}

Integrated nested Laplace approximations (INLA) is a general framework
for performing approximate Bayesian inference in latent Gaussian
models. The INLA
approach is described in 
\cite{rue:martino:chopin:09,Martins:Simpson:Lindgren:Rue:2013}
and has been implemented in the R-INLA package that is publicly available from
the homepage \url{http://www.r-inla.org}.

In fact the INLA approach applies for a wide range of 
latent Gaussian Markov random field models, including
those considered in the present paper, where 
wavelet
coefficients ${\bar{\mathbf w}}$ (the observations) are
conditionally independent given the hidden states ${\bar{\mathbf s}}$
and the unknown
parameters $\v\theta$. 
Here we assume a directional independence GLG model for ${\bar{\mathbf s}}$.
Furthermore, we need to impose a prior distribution on $\v\theta$,
where standard priors have been implemented in R-INLA package that
also allows the user to specify his or her `own prior'.

INLA provides a recipe for fast Bayesian inference using accurate
approximations to $p(\v\theta|{\bar{\mathbf w}})$,
that is, the marginal posterior density 
 for the parameters, and to $p(s_i^{(t)}(\ell)|{\bar{\mathbf w}})$, 
that is, the marginal posterior density 
 for each hidden state $i=1,\ldots,n$ of each tree $t=1,\ldots,k$ and
 each direction $\ell=1,\ldots,d$. 
Thereby we can calculate (approximate) posterior means and
variances as well as other summaries of interest. 
The approximations are based
 on integrated nested Laplace approximations as detailed in 
\cite{rue:martino:chopin:09,Martins:Simpson:Lindgren:Rue:2013}. 
Because the approximations involve
 numerical integration with respect to $\v\theta$, the dimension of
 $\v\theta$ needs to be sufficiently small.

Compared to approximate Bayesian inference using Markov chain Monte Carlo methods, INLA is much faster and accurate, cf.\ \cite{rue:martino:chopin:09}.

\section{Simulation study}
\label{sec:simulations}

To test and compare composite estimates based on the composite EM algorithm with inference based on INLA, we simulate wavelet coefficients from 100 images, each consisting of 256 wavelet trees with three levels, under a homogenous GLG model with parameter values
\begin{equation}
	\label{eq:sim_study_parameters}
	\mu_1 = -1,
	\quad
	\sigma_1 = 1,
	\quad
	\alpha = -1,
	\quad
	\beta = 1,
	\quad
	\kappa = 1.
\end{equation}
This implies that $\mu_2 = -2$, $\mu_3 = -3$, $\sigma_2^2 = 2$ and $\sigma_3^2 = 3$, compare with \textup {(\ref {eq:marginal_moment_relation})}.

Table\nobreakspace \ref {tab:CEM-estimates} shows the average parameter estimates and their standard deviations for estimation with the composite EM algorithm of section Section\nobreakspace \ref {subsec:composite_EM}.
The marginal parameters $\mu$ and $\sigma$ are well estimated, but the transition parameters are not estimated as well.

\begin{table}
	\centering
	\begin{tabular}{ >{$}c<{$} *{3}{>{$}r<{$} @{ } >{$(}l<{$)}} }
		\toprule
		& \multicolumn{2}{c}{$D = 1$} & \multicolumn{2}{c}{$D = 2$} & \multicolumn{2}{c}{$D = 3$}
		\\
		\midrule
		\mu_1    & -0.976 & 0.125 & -0.956 & 0.136 & -0.941 & 0.139
		\\
		\mu_2    & -1.947 & 0.175 & -1.953 & 0.131 & -1.945 & 0.179
		\\
		\mu_3    & -2.956 & 0.270 & -2.968 & 0.209 & -2.952 & 0.277
		\\
		\sigma_1 & 0.937  & 0.291 & 0.915  & 0.284 & 0.891  & 0.275
		\\
		\sigma_2 & 1.843  & 0.314 & 1.847  & 0.323 & 1.813  & 0.343
		\\
		\sigma_3 & 2.823  & 0.429 & 2.824  & 0.435 & 2.813  & 0.460
		\\
		\alpha_1 & -1.562 & 0.332 & -1.662 & 0.310 & -1.654 & 0.360
		\\
		\alpha_2 & -1.707 & 0.340 & -1.748 & 0.323 & -1.716 & 0.325
		\\
		\beta_1  & 0.410  & 0.328 & 0.304  & 0.284 & 0.316  & 0.325
		\\
		\beta_2  & 0.646  & 0.124 & 0.630  & 0.148 & 0.642  & 0.125
		\\
		\kappa_1 & 1.664  & 0.359 & 1.711  & 0.390 & 1.669  & 0.376
		\\
		\kappa_2 & 2.053  & 0.485 & 2.105  & 0.496 & 2.087  & 0.472
		\\
		\bottomrule
	\end{tabular}
	\caption{Results for parameter estimation based on 100 simulations under the homogenous GLG model with parameter values specified in \textup {(\ref {eq:sim_study_parameters})} and obtained using the composite EM algorithm of Section\nobreakspace \ref {subsec:composite_EM}. The $D$ refer to the directional subband. For each parameter the average of 100 parameter estimates is shown together with their standard deviation (show in in parantheses).}
	\label{tab:CEM-estimates}
\end{table}

Table\nobreakspace \ref {tab:INLA_estimates} summarizes the Bayesian inference with INLA.
As a point estimate the marginal posterior mean for each parameter is calculated along with the standard deviation of the 100 estimates, that is,
\begin{gather*}
	\overline{ \E(\theta \,|\, \bar{\v w}) }
	= \frac1{100} \sum_{n=1}^N \E(\theta \,|\, \bar{\v w}_n),
	\\
	\sqrt{\frac1{99} \sum_{n=1}^N \bigl(\E(\theta \,|\, \bar{\v w}_n) - \overline{ \E(\theta \,|\, \bar{\v w}) }\bigr)^2},
\end{gather*}
where $\bar{\v w}_n$ denotes the $n$'th simulated data set.
Note that INLA does not work with variances, but with precisions; that is, we have access to the marginal posterior distributions of $\sigma^{-2}$ and $\kappa^{-2}$.
The estimates in Table\nobreakspace \ref {tab:INLA_estimates} are computed as the average of the inverse posterior means.

\begin{table}
	\centering
	\begin{tabular}{ *{5}{>{$}r<{$} @{ } >{$(}l<{$)}} }
		\toprule
		\multicolumn{2}{c}{$\mu$} & \multicolumn{2}{c}{$\sigma$} & \multicolumn{2}{c}{$\alpha$} & \multicolumn{2}{c}{$\beta$} & \multicolumn{2}{c}{$\kappa$}
		\\
		\midrule
		-0.45 & 0.072 &
		0.010 & 0.051 &
		-1.39 & 0.065 &
		0.82  & 0.052 &
		1.39  & 0.074
		\\
		\bottomrule
	\end{tabular}
	\caption{Results for parameter estimation based on 100 simulations under the homogenous GLG model with parameter values specified in \textup {(\ref {eq:sim_study_parameters})} and obtained using the Bayesian approach of Section\nobreakspace \ref {sec:bayes}. For each parameter the average of 100 parameter estimates is shown together with their standard deviation (show in in parantheses).}
	\label{tab:INLA_estimates}
\end{table}

\begin{figure}
	\centering
	\includegraphics[width=\textwidth]{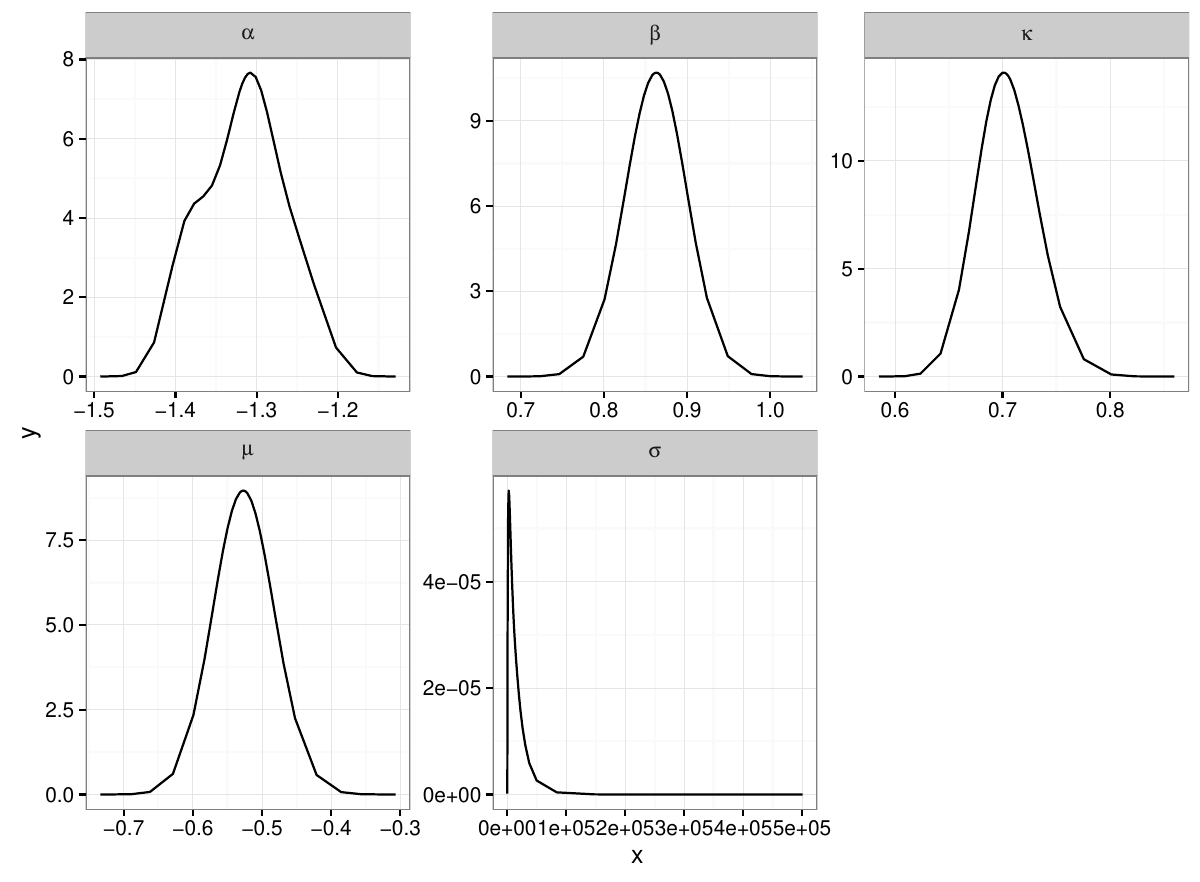}
	\caption{Marginal posterior distributions for each of the parameters in a single simulation. The distributions labellel ``$\sigma$'' and ``$\kappa$'' are actually the posterior distributions for $\sigma^{-2}$ and $\kappa^{-2}$, respectively.}
	\label{fig:sim_study_posterior}
\end{figure}

\begin{figure}
	\centering
	\includegraphics[scale=0.5]{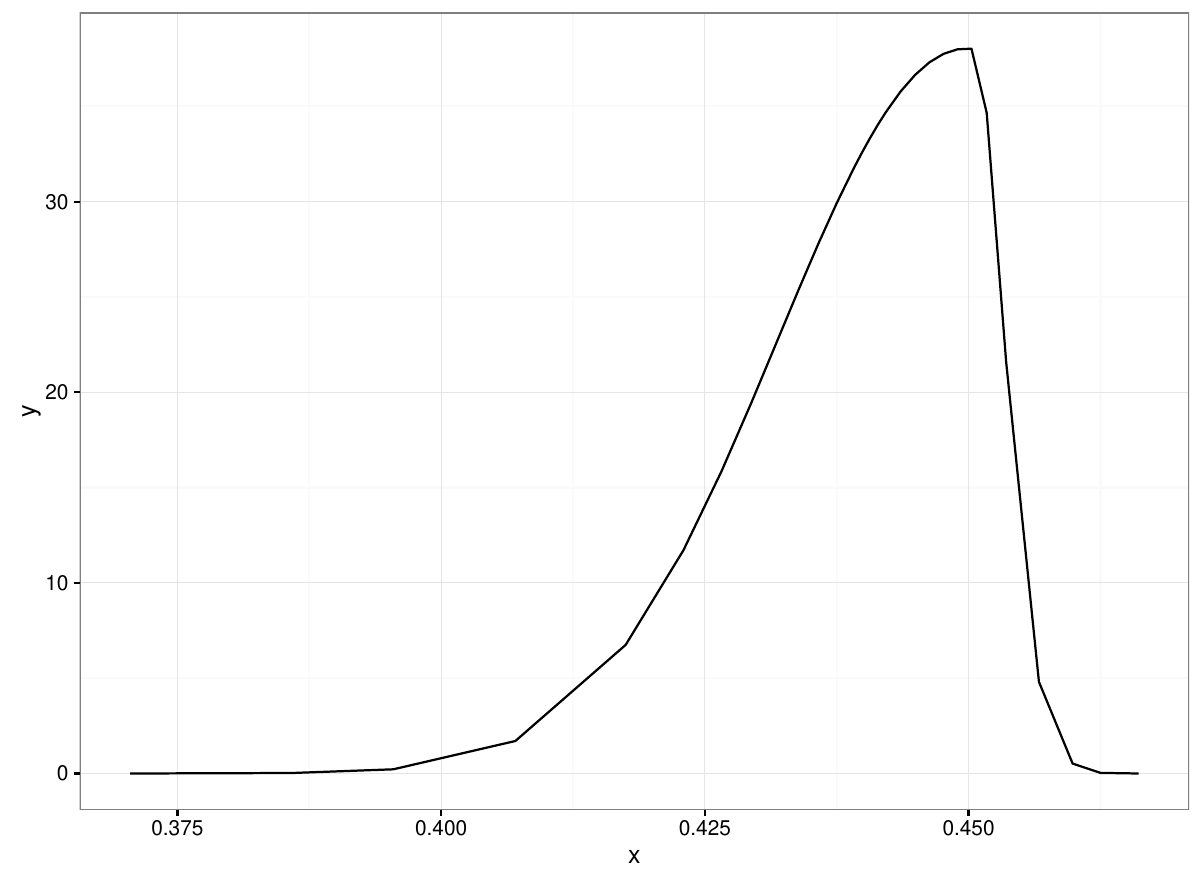}
	\caption{Posterior distribution of $\sigma^{-2}$ from the 'Lena' image in Figure\nobreakspace \ref {fig:test_images}.}
	\label{fig:lena_posterior}
\end{figure}

To illustrate the advantages of the Bayesian approach Figure\nobreakspace \ref {fig:sim_study_posterior} shows the marginal posterior distributions for the parameters of a single simulation.
Unlike the frequentist EM algorithm the Bayesian approach provides information about the uncertainty of the parameters.
One particular thing to notice is the very non-localized distribution for $\sigma^{-2}$.
In our experiments this seems to be a general issue for small data sets -- for a large dataset this posterior distribution is much more localized, see Figure\nobreakspace \ref {fig:lena_posterior} for the posterior distribution of $\sigma^{-2}$ from the 'Lena' image in Figure\nobreakspace \ref {fig:test_images}.

\section{Examples of applications}\label{s:applications}

\begin{figure}[t!]
	\begin{center}
		\includegraphics[width=0.49\textwidth]{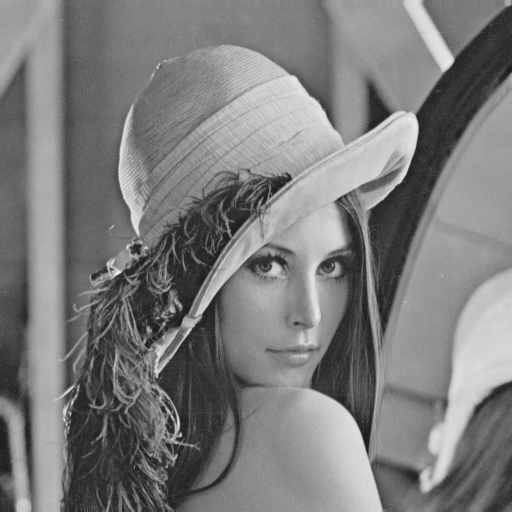}
		\includegraphics[width=0.49\textwidth]{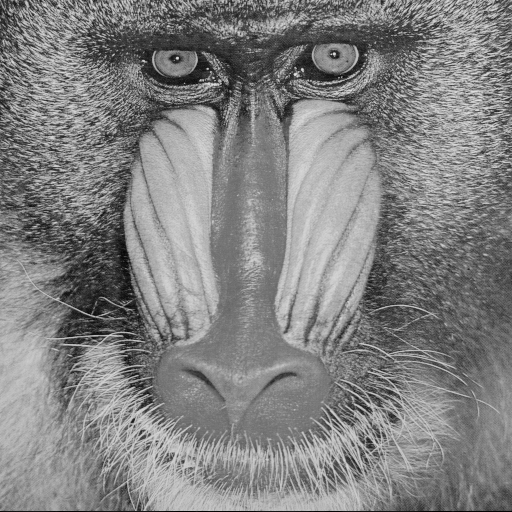}
		\includegraphics[width=0.49\textwidth]{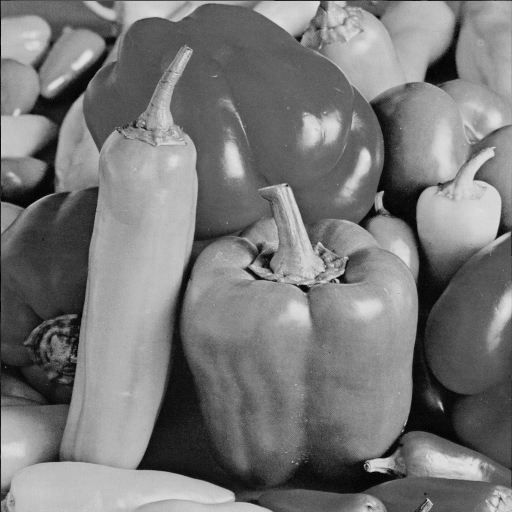}
	\end{center}
	\caption{The three test images: 'Lena', 'mandrill', and 'peppers'.}
	\label{fig:test_images}
\end{figure}

This section demonstrates how the GLG model applies to denoising and edge detection in images.
The examples are meant to illustrate different applications, not to make thorough comparisons with other methods.
We use three test images from the USC-SIPI image database available at \url{http://sipi.usc.edu/database}:
'Lena', 'mandrill', and 'peppers', see Figure\nobreakspace \ref {fig:test_images}.
These images are 512-by-512 pixels represented as 8 bit grayscale with pixel values in the unit interval.

\subsection{Denoising}
\label{s:denoising}

\begin{figure}[t!]
	\begin{center}
		\includegraphics{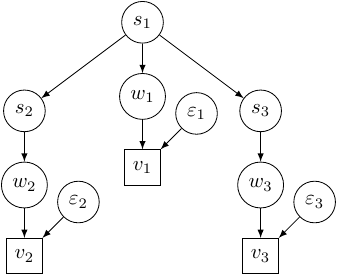}
	\end{center}
	\caption{Graphical model of a binary tree with two levels and noisy observations. 
	The rectangular nodes are observed variables and the round nodes are unobserved variables.}
	\label{fig:noisy_tree_graph}
\end{figure}

Consider an image corrupted with additive white noise, that is, we add 
independent terms to the pixel values from the same zero-mean normal
distribution. Then the procedure for denoising with orthonormal
wavelets works as follows
\begin{equation*}
	\text{noisy data}
	\to
	\text{noisy wavelets}
	\to
	\text{noise-free wavelets}
	\to
	\text{noise-free data}
\end{equation*}
where the distribution and the independence properties of the noise are preserved by the wavelet transform, so the problem of estimating noise-free data boils down to considering each wavelet tree observed with additive white noise:
\begin{equation}\label{e:stjerne}
	v_i^{(t)}(\ell) = w_i^{(t)}(\ell) + \varepsilon_i^{(t)}(\ell),
        \quad \ell=1,\ldots,d,\ i=1,\dots,n,\ t=1,\ldots,k,
\end{equation}
where 
the $\varepsilon_i^{(t)}(\ell)$ are mutually independent, identically $N(0,
\sigma_\varepsilon^2)$-distributed and independent of 
$({\bar{\mathbf w}},{\bar{\mathbf s}})$. 
The dependence structure in a
tree with noisy observations is illustrated in
Figure\nobreakspace \ref {fig:noisy_tree_graph} 
(suppressing the dependence of the indices $t$ and $\ell$).

Below we discuss estimation of each wavelet coefficient $w_i^{(t)}(\ell)$, assuming that the noise variance $\sigma_\epsilon^2$ is known. 
From \textup {(\ref {eq:obs_model})} and \textup {(\ref {e:stjerne})} we obtain the conditional density
\begin{equation*}
	p(\mathbf w | \mathbf v, \v s, \v\theta) = \prod_{i=1}^n \prod_{\ell=1}^d p\bigl(w_i(\ell) | v_i(\ell), s_i(\ell), \v\theta\bigr).
\end{equation*}
In the sequel, to simplify the notation, 
we drop the indices and denote $w_i^{(t)}(\ell)$, $v_i^{(t)}(\ell)$, 
$s_i^{(t)}(\ell)$ by $w$, $v$, $s$, respectively, where $s\sim
N(\mu,\sigma^2)$, cf.\ \eqref{eq:marginal_node_pdf}. 

In the frequentist setup, we estimate $w$ by the conditional mean
\begin{equation}
	\label{eq:freq_post_mean}
	\E[w | v, \v\theta] 
	= \frac{v}{c(v | \mu, \sigma^2)}
	\int \frac{ \exp(s) }{(\exp(s) + \sigma_\varepsilon^2)^{3/2}} 
	\exp\Biggl(-\frac{1}{2}\Biggl[\frac{v^2}{\exp(s) + \sigma_\varepsilon^2} + \frac{(s - \mu)^2}{\sigma^2}\Biggr]\Biggr) \diff s
\end{equation}
as derived later in \eqref{AHA} and
where $c(v | \mu, \sigma^2)$ is defined in
Appendix~\ref{app:denoise}. Here 
we use the Gauss-Hermite quadrature rule for approximating the integral. 
Furthermore, since \textup {(\ref {eq:freq_post_mean})} depends only on the
parameters $\mu$ and $\sigma^2$, we
replace these 
by the estimates obtained by the composite EM algorithm. 

In the Bayesian setup, we 
work with the posterior density $p(w | v)$ from which we can calculate various point estimates.
We have 
\begin{align*}
p(w | v)= \int p(w | v, s) p(s | v)\d s,\quad 
\E(w | v)= \int \E(w | v, s) p(s | v)\d s,
\end{align*}
where $p(s | v)$ is calculated using the INLA approach.  
For fixed values of $v$ and $s$, we have $p(w | v, s)
	\propto p(w | s) p(v | w)$, where
$p(w | s) \sim N(0, \exp(s))$ and $p(v | w) \sim N(w,
\sigma_\epsilon^2)$. Therefore
\begin{equation*}
	p(w | v, s)
	\sim N\biggl( \frac{v \exp(s)}{\sigma_\epsilon^2 + \exp(s)}, \frac{\sigma_\epsilon^2 \exp(s)}{\sigma_\epsilon^2 + \exp(s)} \biggr)
\end{equation*}
and so we can evaluate e.g.
\begin{equation}\label{eq:denoising_posterior_pdf}
	\E(w | v) 
	= v \int \frac{\exp(s)}{\sigma_\epsilon^2 + \exp(s)} p(s | v) \d s
\end{equation}
by numerical integration. 

We applied the two denoising schemes with a 3 level wavelet transform using the Daubechies 4 filter to noisy versions of the three test images in Figure\nobreakspace \ref {fig:test_images}. 
To estimate the performance of a denoising scheme, we calculated the peak signal-to-noise ratio (PSNR) in decibels between a test image $I$ and a noisy or cleaned image $J$.
For images of size $N\times N$, the PSNR in decibels is defined by
\begin{equation*}
	\label{eq:PSNR_def}
	\text{PSNR} = 20 \log_{10} \frac{N (\max\{I(x)\} - \min\{I(x)\})}{\|I - J\|}
\end{equation*}
where the maximum and the minimum are over all pixels $x$ and where 
$\|\cdot\|$ is the Frobenius norm.
Table\nobreakspace \ref {tab:denoising_results} shows for the test images and different
noise levels $\sigma_\varepsilon$, the PSNR between each test image
and its noisy or denoised version obtained by the frequentist approach
under the directional independence GLG model or by the Bayesian
approach under the homogeneous GLG model.
The Bayesian results yields the lowest PSNR values, but they are also
based on a more parsimonious model. Moreover,  
the posterior mean \eqref{eq:denoising_posterior_pdf}
and the posterior median based on $p(w | v)$ were almost identical.
\begin{table}[t!]
	\centering
	\caption{For the three test images and three noise levels, peak signal-to-noise ratios in dB between the image and its noisy version ot its denoised version obtained using either the Gaussian Finite Mixture model and the EM algorithm from \cite{crouse:nowak:baraniuk:98}, the GLG model and the composite EM algorithm, or the homogeneous GLG model and INLA. 
In the latter case, the PSNR is calculated using the median of the posterior image.
	For each image, a three level Daubechies 4 wavelet transform is used.}
	\label{tab:denoising_results}
	\begin{tabular}{*{5}{c}}
		\toprule
		& & \multicolumn{3}{c}{PSNR} \\
		\cmidrule{3-5}
		test image & \smash{\shortstack{noise\\ level $\sigma_\varepsilon$}} & noisy & Freq. & Bayes 
		\\
		\midrule
				 & 0.10 & 18.76 & 27.93 & 23.57 \\
		Lena 	 & 0.15 & 15.44 & 26.18 & 20.68 \\
				 & 0.20 & 13.17 & 24.72 & 18.77 \\
		\midrule
				 & 0.10 & 19.18 & 23.39 & 22.47 \\
		Mandrill & 0.15 & 15.77 & 21.61 & 19.70 \\
				 & 0.20 & 13.49 & 20.52 & 18.09 \\
		\midrule
				 & 0.10 & 19.18 & 27.96 & 24.00 \\
		Peppers  & 0.15 & 15.83 & 25.87 & 21.08 \\
				 & 0.20 & 13.57 & 24.41 & 19.18 \\
		\bottomrule
	\end{tabular}
\end{table}

\begin{figure}[t!]
	\begin{center}
		\includegraphics[width=\textwidth]{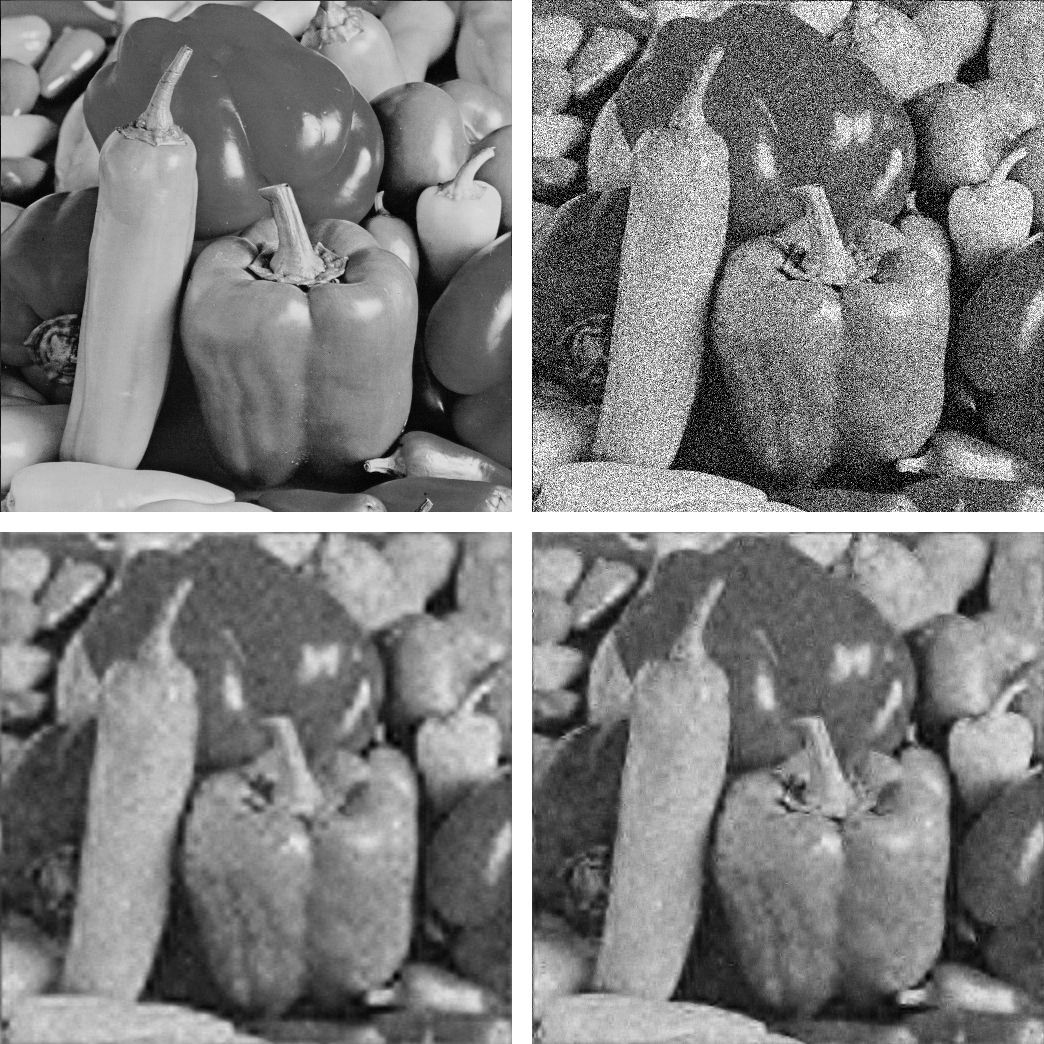}
	\end{center}
	\caption{Denoising results for the peppers image from Table\nobreakspace \ref {tab:denoising_results} when the standard deviation of the noise is 0.20.  
	Top left panel: The original image.  
	Top right panel: The noisy image (PSNR is 13.57).  
	Bottom image: The noisy image cleaned using the GLG model and the composite EM algorithm (PSNR is 24.41).}
	\label{fig:denoising_example}
\end{figure}

An example of the visual appearance of denoising using frequentist means is seen in Figure\nobreakspace \ref {fig:denoising_example}. 

\begin{figure}[t!]
  \begin{center}
	  \includegraphics[width=\textwidth]{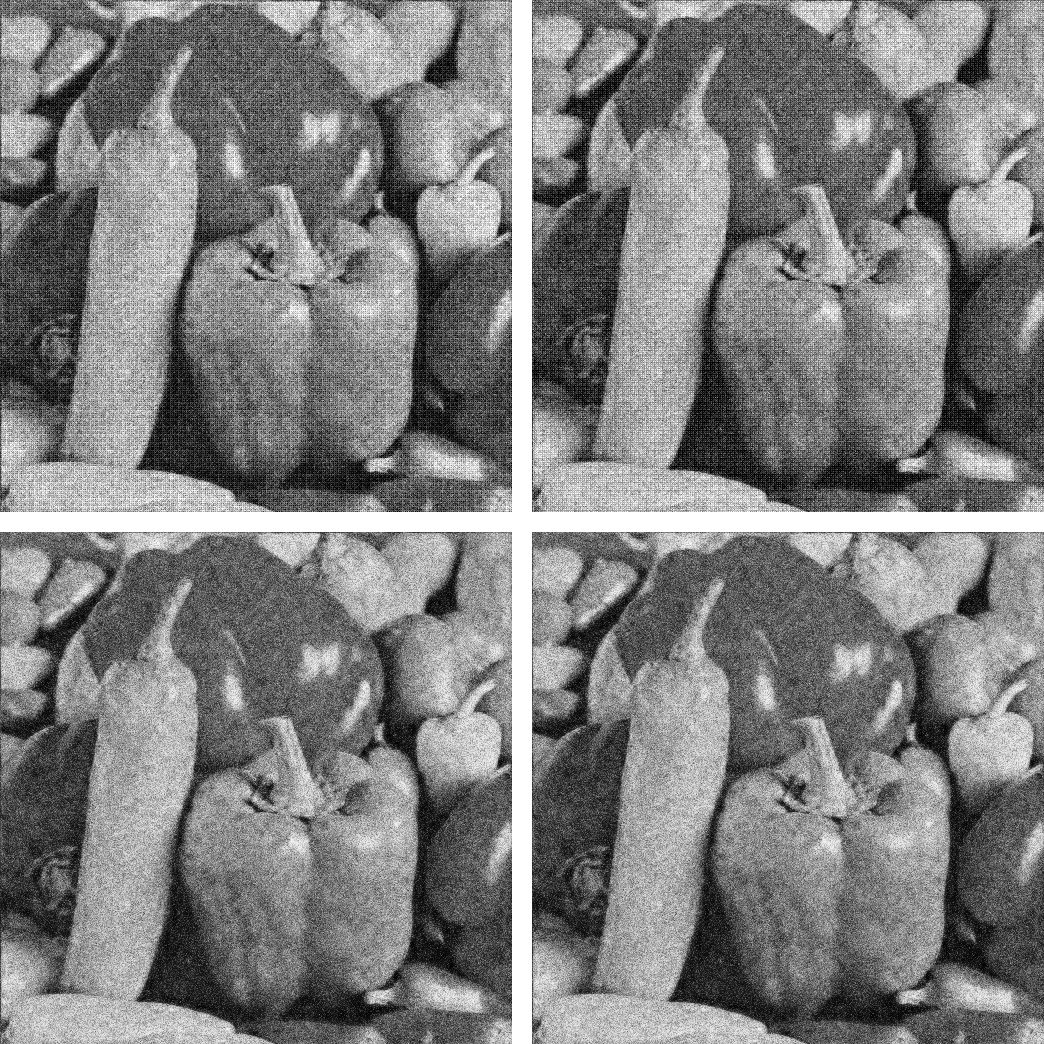}
  \end{center}
  \caption{Denoising the 'peppers' image using the posterior distribution \textup {(\ref {eq:denoising_posterior_pdf})} and INLA.
  	The original and noisy images are seen in Figure\nobreakspace \ref {fig:denoising_example}.
	The top left, top right, and bottom left images are based on the 25\%, 75\%, and 50\% quantiles of the posterior distribution, respectively (the PSNRs are 16.38, 16.41, and 19.18, respectively).
	The bottom right image is based on the mean of the posterior distribution (PSNR is 19.15).
	The posterior mean and median are almost identical.}
  \label{fig:Bayesian_denoising_example}
\end{figure}

The median (the 50\% quantile) of the posterior distribution is only one possible point estimate of the posterior distribution. 
However, using other quantiles or the posterior mean are not providing better results, see Figure\nobreakspace \ref {fig:Bayesian_denoising_example}. 

\subsection{Edge detection}
\label{sec:edge_detection}

Edge detection in an image is performed by labelling each pixel as being either an edge or a non-edge. 
Using the wavelet transform for this task has the advantage that wavelet coefficients are large near edges and small in the homogeneous parts of an image;
the difficulty lies in quantifying `large' and `small'.
Another advantage is that a multiresolution analysis allows us to
search for edges that are present at only selected scales of the
image, that is, 
edges that are neither too coarse nor too fine.

This section 
discusses how to
label the wavelet coefficient $w_i^{(t)}(\ell)$ by an indicator
variable $f_i^{(t)}(\ell)$ so that $f_i^{(t)}(\ell) = 1$ indicates that
$w_i^{(t)}(\ell)$ is `large' and $f_i^{(t)}(\ell) = 0$ indicates that $w_i(\ell)^{(t)}$ is `small'.
Our procedure is inspired by that of \cite{Sun:Gu:Chen:Zhang:2004} 
which uses the 2-state Gaussian Finite Mixture (GFM) model of
\cite{crouse:nowak:baraniuk:98}. Below
we recap this labelling algorithm, modify it to the case of our GLG
model, and discuss how to transfer the labels $f_i^{(t)}(\ell)$ to the pixels.
Finally, we show examples of both the original procedure using the GFM
model and our modified procedure using the GLG model.

For our brief description of the edge detection algorithm in
\cite{Sun:Gu:Chen:Zhang:2004} it is sufficient to consider the 2-state
GFM model based on a conditional independence structure as in
\textup {(\ref {eq:full_model})}, but where the hidden states take only binary values, indicating whether the associated wavelet coefficient is large or not.
The labelling in \cite{Sun:Gu:Chen:Zhang:2004} consists of three steps.
First, using the EM algorithm 
an estimate ${\widehat{\v\theta}}$ of the parameter vector $\v\theta$ of a 2-state GFM model is obtained from the data ${\bar{\mathbf w}}$.
Second, using an empirical Bayesian approach, the maximum a posteriori (MAP) estimate of the hidden states
\begin{equation}
	\label{eq:viterbi_eq}
	\widehat{\mathbf s^{(t)}} = \argmax_{\mathbf s^{(t)}} p(\mathbf s^{(t)} |
        \mathbf w^{(t)}, {\widehat{\v\theta}}) = \argmax_{\mathbf s^{(t)}}
        p(\mathbf s^{(t)}, \mathbf w^{(t)} | {\widehat{\v\theta}}),\quad t = 1,\ldots,k,
\end{equation}
 is computed.
Finally, they define $f_i^{(t)}(\ell) = {\widehat{s^{(t)}_i(\ell)}}$.

The idea of labelling wavelet coefficients with the GLG model is
overall the same as presented above for the GFM model, with the
differences arising from the continuous nature of the hidden states
and from 
different algorithms being applied for parameter 
estimation and state estimation.
First, an estimate ${\widehat{\v\theta}}$ of the parameter vector of the GLG model is obtained.
Second, in analogy with \textup {(\ref {eq:viterbi_eq})} we compute the MAP estimate ${\widehat{\mathbf s^{(t)}}}$ by noting that
\begin{equation}
	\label{eq:GLG_state_estimation_factorization}
	p(\mathbf s^{(t)}, \mathbf w^{(t)}|{\widehat{\v\theta}}) 
	= p(\mathbf s^{(t)}|{\widehat{\v\theta}}) \prod_{i=1}^n \prod_{\ell=1}^d p\bigl(w_i^{(t)}(\ell) | s^{(t)}_i(\ell)\bigr)
\end{equation}
where $p(\mathbf s^{(t)}|{\widehat{\v\theta}})$ is a multidimensional Gaussian density
function those estimated mean vector and precision matrix are denoted 
${\widehat{\bm\mu}}$ and 
$\widehat\Delta$, respectively.
The log of \textup {(\ref {eq:GLG_state_estimation_factorization})} and its gradient vector and Hessian matrix with respect to $\mathbf s^{(t)}$ are
\begin{gather*}
	\log p(\mathbf s^{(t)}, \mathbf w^{(t)}|{\widehat{\v\theta}}) \equiv -\frac12 \Bigl\{(\mathbf s^{(t)} - {\widehat{\bm\mu}})^\top \widehat\Delta (\mathbf s^{(t)} - {\widehat{\bm\mu}}) + \sum_{i=1}^n \sum_{\ell=1}^d \bigl(w^{(t)}_i(\ell)^2 \exp(-s^{(t)}_i(\ell)) + s^{(t)}_i(\ell)\bigr)\Bigr\},
	\\
	\nabla\log p(\mathbf s^{(t)}, \mathbf
        w^{(t)}|{\widehat{\v\theta}}) = -\widehat\Delta (\mathbf
        s^{(t)} - {\widehat{\bm\mu}}) + \frac12
        \bigl[w^{(t)}_i(\ell)^2 \exp(-s^{(t)}_i(\ell)) -
        1\bigr]_{\ell=1,\ldots,d,\, i=1,\ldots, n},
	\\
	H\bigl(\log p(\mathbf s^{(t)}, \mathbf
        w^{(t)})\bigr|{\widehat{\v\theta}}) = -\widehat\Delta -
        \frac12 \diag\bigl( w^{(t)}_i(\ell)^2 \exp(-s^{(t)}_i(\ell)),
        \, \ell=1,\ldots,d,\, i=1,\ldots, n \bigr),
\end{gather*}
where $\equiv$ means that an additive term which is not depending on $\mathbf s^{(t)}$ has been omitted in the right hand side expression.
The Hessian matrix is strictly negative definite for all $(\mathbf s^{(t)}, \mathbf w^{(t)})$ with $\mathbf w^{(t)} \neq \v 0$, and hence $\widehat{\mathbf s^{(t)}}$ can be found by solving $\nabla\log p(\mathbf s^{(t)}, \mathbf w^{(t)}|{\widehat{\v\theta}}) = 0$ using standard numerical tools.
Finally, observe that if the estimate ${\widehat{s^{(t)}_i(\ell)}}$ is
deemed to be large in the estimated distribution $N({\widehat{\mu_{\rho(i)}}}, {\widehat{\sigma^2_{\rho(i)}}})$ for $s^{(t)}_i(\ell)$, then we expect $w^{(t)}_i(\ell)$ to be `large'.
Therefore, denoting $z_p$ the $p$-quantile in
$N({\widehat{\mu_{\rho(i)}}}, {\widehat{\sigma_{\rho(i)}^2}})$, with e.g.\ $p = 0.9$, we define $f^{(t)}_i(\ell) = 1$ if ${\widehat{s^{(t)}_i(\ell)}}\geq z_p$ and $f^{(t)}_i(\ell) = 0$ otherwise.

It remains to clarify how to transfer $f_i^{(t)}(\ell)$ (defined by
one of the two methods above) to the pixel domain (this issue is not
discussed in \cite{Sun:Gu:Chen:Zhang:2004}). That is, we want specify
a binary iamge with pixel values $e_j$ indicating whether pixel $j$ is
part of an edge or not. Applying the inverse wavelet transform (with
low-pass coefficients that are zero) to the label trees $\mathbf
f^{(t)}$, $t=1,\ldots,k$, give an image with pixel values $\tilde
e_j$, say. Since the wavelet transform 
does not necessarily map binary values to binary values, we define 
\begin{equation*}
	e_j =
	\begin{dcases}
		1 & \text{if } \widetilde e_j \neq 0,\\
		0 & \text{otherwise}.
	\end{dcases}
\end{equation*}
The $e_j$'s are sensitive to the choice of wavelet transform, and
using e.g.\ the Haar wavelet results in thin edges.

As mentioned, the multiresolution analysis of the wavelet transform
allows us to consider edges with properties that are present at only specific scales.
To exclude edges at a given level $r$, direction $\ell$ and tree $t$ in the wavelet transform, we simply
modify $f^{(t)}_i(\ell)$ by setting $f^{(t)}_i(\ell) = 0$ whenever $\l(i) = r$.

\begin{figure}[t!]
	\begin{center}
		\includegraphics[width=\textwidth]{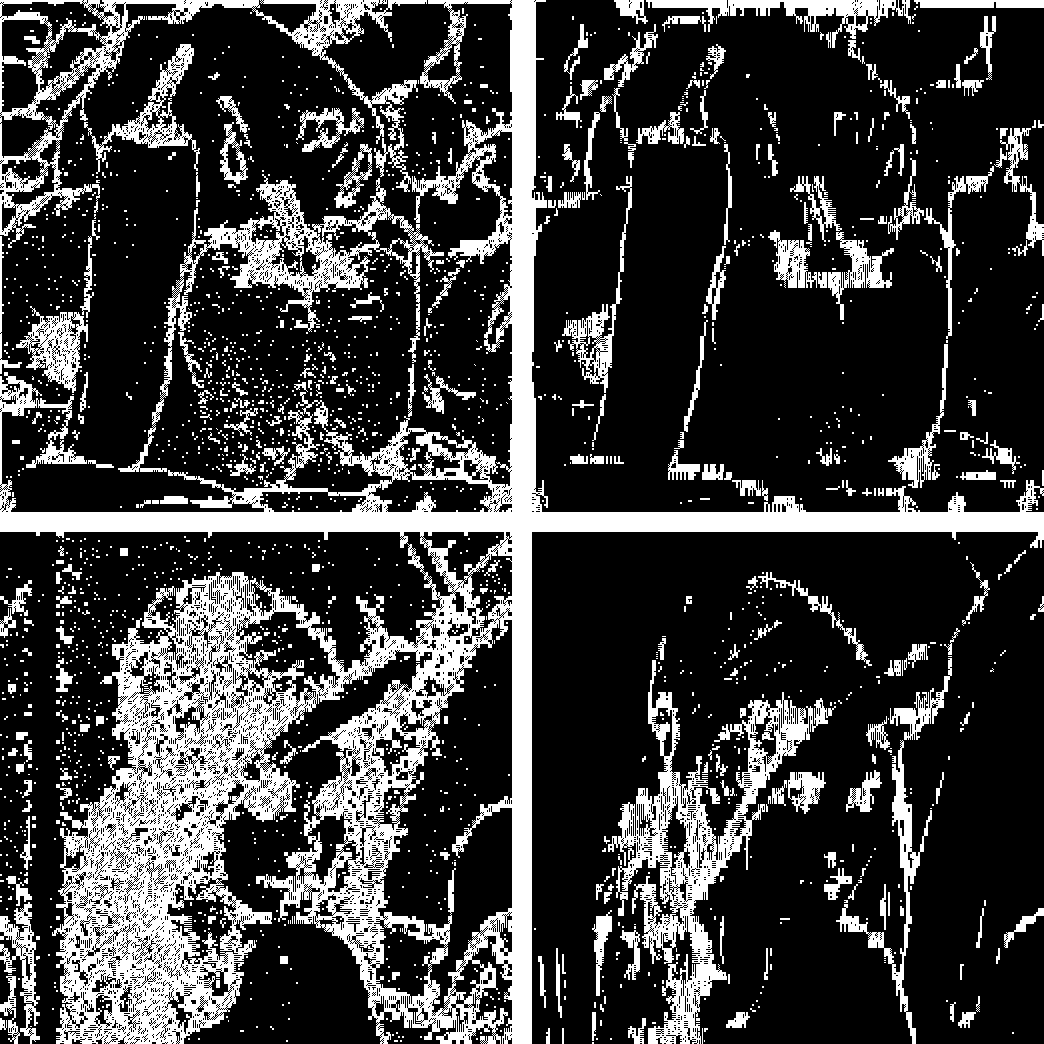}
	\end{center}
	\caption{Examples of edge detection of the 'Lena' and
          'peppers' images using the method from
          \protect\cite{Sun:Gu:Chen:Zhang:2004} (left column) and our
          variant that uses the GLG model (right column), when using
	a 3 level Haar wavelet transform and the $90\%$-quantile 
for thresholding with the GLG model, and  considering
only the finest level of the wavelet transform.}
	\label{fig:edge_detection}
\end{figure}

Figure\nobreakspace \ref {fig:edge_detection} compares the results of the two edge
detection algorithms, when we only used the finest scale in the
wavelet transform, and where in the case of our method the directional
independence GLG model and the composite EM algorithm were used.
Our method classified fewer pixels as edges, and comparing with Figure\nobreakspace \ref {fig:test_images} it appears that the GFM-based procedure includes many superfluous pixels.


\section*{Acknowledgment}

This research is supported by the Danish Council for Independent Research | Natural Sciences, grant 12-124675, "Mathematical and Statistical Analysis of Spatial Data", and by the Centre for Stochastic Geometry and Advanced Bioimaging, funded by a grant (8721) from the Villum Foundation. 
We are grateful to H{\aa}vard Rue for help with INLA. 
We thank Peter Craigmile, Morten Nielsen, and Mohammad Emtiyaz Khan for helpful discussions.
We also thank the reviewers for their constructive feedback on earlier versions of the manuscript that improved our contribution significantly.

\appendix

\section{Moments and estimating equations}
\label{app:moments}

This appendix specifies the moment estimates briefly discussed in Section\nobreakspace \ref {sec:moment_relation}.

Let $i$ be a node on level $r=r(i)$.
Using \textup {(\ref {eq:obs_model})} and\nobreakspace  \textup {(\ref {eq:marginal_node_pdf})}, conditioning on the hidden states and exploiting the conditional independence structure, we obtain
\begin{align}
	\eta_{r,\ell}^{(2)} 
	& := \E\bigl[w_i(\ell)^2\bigr] 
	= \exp\Bigl(\mu_{r,\ell} + \sigma_{r,\ell,\ell}/2\Bigr), 
	\label{eq:eta2}
	\\
	\eta_{r,\ell}^{(4)}
	& := \E\bigl[w_i(\ell)^4\bigr] 
	= 3\exp\Bigl(2\mu_{r,\ell} + 2\sigma_{r,\ell,\ell}\Bigr), 
	\label{eq:eta4}
\end{align}
for $1\leq\ell \leq d$, 
\begin{align}
	\eta_{r,\ell, \ell'}^{(2,2)} 
	& := \E\bigl[w_i(\ell)^2 w_i(\ell')^2\bigr]
	\nonumber
	\\
	& = \exp\Bigl(\mu_{r,\ell} + \mu_{r,\ell'} + \sigma_{r,\ell, \ell}/2 + \sigma_{r,\ell',\ell'}/2 + \sigma_{r,\ell,\ell'}\Bigr), 
	\label{eq:eta22}
\end{align}
for $1\leq\ell < \ell' \leq d$, and
if $r < L$ and $j \in c(i)$, then
\begin{align}
	\xi_{r,\ell,\ell'}^{(2,2)}
	& := \E\bigl[w_i(\ell)^2 w_j(\ell')^2\bigr]
	\nonumber
	\\
	& = \exp\bigl( \mu_{r,\ell} + \mu_{r+1,\ell'} + \sigma_{r,\ell,\ell}/2 + \sigma_{r+1,\ell',\ell'}/2 + \sum_{k=1}^d \beta_{r,\ell,k} \sigma_{r,k,\ell'} \bigr),
	\label{eq:xi22}
\end{align}
for $1\leq \ell, \ell' \leq d$. Let $n_r$ denote the number of nodes
on level $r$. Using the unbiased estimates of
$\eta_{r,\ell}^{(2)}$, $\eta_{r,\ell}^{(4)}$, 
$\eta_{r,\ell,\ell'}^{(2,2)}$ given by
\begin{align*}
	\eta_{r,\ell}^{(a)}
	& = \frac1{k n_r} \sum_{t=1}^k \sum_{i: \l(i) = r} \sum_{j\in c(i)} \bigl(w_j(\ell)^{(t)}\bigr)^a,
	\quad a = 2,4,
	\\
	\eta_{r, \ell, \ell'}^{(2,2)}
	& = \frac1{k n_r} \sum_{t=1}^k \sum_{i: \l(i) = r} \sum_{j\in c(i)} \bigl(w_j(\ell)^{(t)}\bigr)^2 \bigl(w_j(\ell')^{(t)}\bigr)^2,
	\\
	\xi_{r, \ell, \ell'}^{(2,2)}
	& = \frac1{k n_r} \sum_{t=1}^k \sum_{i: \l(i) = r} \sum_{j\in c(i)} \bigl(w_i(\ell)^{(t)}\bigr)^2 \bigl(w_j(\ell')^{(t)}\bigr)^2,
\end{align*} 
 \textup {(\ref {eq:eta2})}- \textup {(\ref {eq:eta22})}  can be solved to obtain estimates 
${\widehat{\v\mu_r}}$ and ${\widehat{\m\Sigma_r}}$. Combining
the estimates $\widehat{\v\mu_r}$, $\widehat{\m\Sigma_r}$,
$\widehat{\xi}_{r,\ell,\ell'}^{(2,2)}$ with \textup {(\ref {eq:xi22})} we obtain
an estimate ${\widehat{\m\Beta_r}}$. Finally, combining these estiamtes
with \textup {(\ref {eq:marginal_moment_relation})} we obtain
estimates ${\widehat{\v\alpha_r}}$ and ${\widehat{\m\Kappa_r}}$. 

The estimating equations  \textup {(\ref {eq:eta2})}- \textup {(\ref {eq:xi22})}  do not guarantee that the estimated covariance matrices are strictly positive definite.
Therefore the estimates of $\m\Sigma_r$ and $\m\Kappa_r$ are replaced by the positive definite matrices that are closest in the Frobenius norm.
This choice of norm is due to the simple analytical solution of the positive definite approximant \cite{Higham:1988}.

\section{EM algorithm for marginal likelihoods}
\label{app:compEM}

The EM algorithm \cite{dempster:laird:rubin:77,Gao:Song:2011} is an iterative estimation procedure which applies for steps \ref {item:root_EM} and\nobreakspace  \ref {item:transfer_EM} in Section\nobreakspace \ref {subsec:composite_EM} as described below.

We start by noticing that the conditional density of $\v s_1$ given $\v w_1$ is
\begin{align}
	\MoveEqLeft 
	p(\v s_1 | \v w_1, \v\mu_1, \m\Sigma_1) 
	= \frac{ p(\v s_1, \v w_1 \,|\, \v\mu_1, \m\Sigma_1) }{ q(\v w_1 \,|\, \v\mu_1, \m\Sigma_1) }
	\nonumber
	\\
	& \propto \exp\biggl(-\frac12 \biggl\{\biggl[\sum_{\ell=1}^d \frac{w_1(\ell)^2}{\exp(s_1(\ell))} + s_1(\ell)\biggr] + (\v s_1 - \v\mu_1)^\T \m\Sigma_1^{-1} (\v s_1 - \v\mu_1)\biggr\}\biggr)
	\label{eq:conditional_pdf_root}
\end{align}
where in the expression on the right hand side we have omitted a factor which does not depend on the argument $\v s_1$ of the conditional density.
Note also that for $\l(i) = r < L$, the conditional density of $\v s_{i,c(i)}$ given $\v w_{i,c(i)}$ is
\begin{align}
	\MoveEqLeft
	p(\v s_{i,c(i)} \,|\, \v w_{i,c(i)}, \v\mu_{r-1}, \m\Sigma_{r-1}, \v\alpha_r, \m\Beta_r, \m\Kappa_r) 
	= \frac{ p(\v s_{i,c(i)}, \v w_{i,c(i)} \,|\, \v\mu_{r-1}, \m\Sigma_{r-1}, \v\alpha_r, \m\Beta_r, \m\Kappa_r) }
	{ q(\v w_{i,c(i)} \,|\, \v\mu_{r-1}, \m\Sigma_{r-1}, \v\alpha_r, \m\Beta_r, \m\Kappa_r) }
	\nonumber
	\\
	& \propto \exp\biggl(-\frac12 \biggl\{\biggl[\sum_{\ell=1}^d \frac{w_i(\ell)^2}{\exp(s_i(\ell))} + s_i(\ell)\biggr] + (\v s_i - \v\mu_{r-1})^\T \m\Sigma_{r-1}^{-1} (\v s_i - \v\mu_{r-1})
	\nonumber
	\\
	& + \sum_{j\in c(i)} \biggl[\sum_{\ell=1}^d \frac{w_j(\ell)^2}{\exp(s_j(\ell))} + s_j(\ell)\biggr] + (\v s_j - \v\alpha_r - \m\Beta_r \v s_i)^\T \m\Kappa_r^{-1} (\v s_j - \v\alpha_r - \m\Beta_r \v s_i)\biggr\}\biggr)
	\label{eq:conditional_pdf_transfer}
\end{align}

In step \ref {item:root_EM}, suppose $(\v\mu_1', \m\Sigma_1')$ is a current estimate.
We consider the conditional expectation with respect to \textup {(\ref {eq:conditional_pdf_root})} when $(\v\mu_1, \m\Sigma_1)$ is replaced by $(\v\mu_1', \m\Sigma_1')$.
Then the next estimate for $(\v\mu_1, \m\Sigma_1)$ is the maximum point for the conditional expectation of the log-likelihood which is based on both $\overline{\v w}_1$ and $\overline{\v s}_1$;
this log-likelihood is given by
\begin{equation*}
	\sum_{t=1}^k \log p(\v w_1^{(t)}, \v s_1^{(t)} \,\big|\, \v\mu_1, \m\Sigma_1)
	\equiv -\frac12 \sum_{t=1}^k \Bigl[\log\det(\m\Sigma_1) + (\v s_1 - \v\mu_1)^\T \m\Sigma_1^{-1} (\v s_1 - \v\mu_1)\Bigr]
\end{equation*}
where $\equiv$ means that an additive term which is not depending on $(\v\mu_1, \m\Sigma_1)$ has been omitted in the right hand side expression, cf.\ \textup {(\ref {eq:jointpdf_root})}.
It is well-known from the theory of estimation in the multivariate Gaussian distribution that the maximum point is given by
\begin{gather}
	\widehat{\v\mu}_1 
	= \frac1k \sum_{t=1}^k \E\Bigl[\v s_1^{(t)} \,|\, \v w_1^{(t)}, \v\mu_1', \m\Sigma_1'\Bigr],
	\label{eq:CEM_mu0}
	\\
	\widehat{\v\Sigma}_1 
	= \biggl\{\frac1k \sum_{t=1}^k \E\Bigl[\v s_1^{(t)} \bigl(\v s_1^{(t)}\bigr)^\T \,|\, \v w_1^{(t)}, \v\mu_1', \m\Sigma_1'\Bigr]\biggr\} - \widehat{\v\mu}_1 \widehat{\m\mu}_1^\T,
	\label{eq:CEM_sigma0}
\end{gather}
where the conditional expectations are calculated using \textup {(\ref {eq:conditional_pdf_root})}.
As mentioned in Section\nobreakspace \ref {sec:marginal_likelihoods}, the integrals are calculated using the Gauss-Hermite quadrature rule.
The iteration is repeated with $(\v\mu_1', \m\Sigma_1') = (\widehat{\v\mu}_1, \widehat{\m\Sigma}_1)$ until convergence is effectively obtained, whereby a final estimate $(\widehat{\v\mu}_1, \widehat{\m\Sigma}_1)$ is returned.

In step \ref {item:transfer_EM}, suppose $(\v\alpha_r', \m\Beta_r', \m\Kappa_r')$ is a current estimate, which we use together with the estimate $(\widehat{\v\mu}_{r}, \widehat{\m\Sigma}_{r})$ to obtain the next estimate for $(\v\alpha_r, \m\Beta_r, \m\Kappa_r)$:
Replacing $(\v\mu_{r}, \m\Sigma_{r}, \v\alpha_r, \m\Beta_r, \m\Kappa_r)$ by $(\widehat{\v\mu}_{r}, \widehat{\m\Sigma}_{r}, \v\alpha_r', \m\Beta_r', \m\Kappa_r')$, this estimate is the maximum point for the conditional expectation with respect to \textup {(\ref {eq:conditional_pdf_transfer})} of each term in the sum
\begin{align*}
	\MoveEqLeft
	\sum_{t=1}^k \sum_{i: \l(i) = r} \log p(\v w_{i,c(i)}^{(t)}, \v s_{i,c(i)}^{(t)} \,\big|\, \widehat{\v\mu}_{r}, \widehat{\m\Sigma}_{r}, \v\alpha_r, \m\Beta_r, \m\Kappa_r)
	\\
	& \equiv -\frac12 \sum_{t=1}^k \sum_{i: \l(i) = r} \sum_{j\in c(i)} \Bigl[ \log\det(\m\Kappa_r) + (\v s_j - \v\alpha_r - \m\Beta_r \v s_i)^\T \m\Kappa_r^{-1} (\v s_j - \v\alpha_r - \m\Beta_r \v s_i) \Bigr],
\end{align*}
where additive terms which do not depend on $(\v\alpha_r, \m\Beta_r, \m\Kappa_r)$ have been omitted, cf.\ \textup {(\ref {eq:jointpdf_transfer})}.
For $1\leq r < L$, let $n_r = \#\{i \,:\, \l(i) = r\}$ and
\begin{equation*}
	\overline{\v s}(r) 
	= \frac1{k n_r} \sum_{t=1}^k \sum_{i : \l(i)=r} \E\Bigl[ \v s_k^{(t)} \,\Big|\, \v w_{i,c(i)}^{(t)}, \widehat{\v\mu}_{r}, \widehat{\m\Sigma}_{r}, \v\alpha_r', \m\Beta_r', \m\Kappa_r' \Bigr].
\end{equation*}
Then, by the theory of multiple linear regression, the maximum point is given by
\begin{align}
	\widehat{\m\Beta}_r & =
	\biggl\{\sum_{t=1}^k \sum_{i: \l(i) = r} |c(i)|
	\E\Bigl[\v s_i^{(t)} \bigl(\v s_i^{(t)}\bigr)^\T - \overline{\v s}(r) \overline{\v s}(r)^\T \,\Big|\, \v w_{i,c(i)}^{(t)}, \widehat{\v\mu}_{r}, \widehat{\m\Sigma}_{r}, \v\alpha_r', \m\Beta_r', \m\Kappa_r'\Bigr]\biggr\}^{-1}
	\nonumber
	\\
	& \quad \sum_{t=1}^k \sum_{i: \l(i) = r} \sum_{j\in c(i)} 
	\E\Bigl[\v s_j^{(t)} \bigl(\v s_i^{(t)} - \overline{\v s}(r)\bigr)^\T \,\Big|\, \v w_{i,c(i)}^{(t)}, \widehat{\v\mu}_{r}, \widehat{\m\Sigma}_{r}, \v\alpha_r', \m\Beta_r', \m\Kappa_r'\Bigr]
	\label{eq:CEM_beta}
	\\
	\widehat{\v\alpha}_r
	& = \overline{\v s}(r+1) - \widehat{\m\Beta}_r \overline{\v s}(r),
	\label{eq:CEM_alpha}
	\\
	\widehat{\m\Kappa}_r
	& = \frac1{k n_r} \sum_{t=1}^k \sum_{i: \l(i) = r} \frac1{|c(i)|} \sum_{j\in c(i)}
	\nonumber
	\\
	&\quad \E\Bigl[ (\v s_j^{(t)} - \widehat{\m\Beta}_r \v s_i^{(t)}) (\v s_j^{(t)} - \widehat{\m\Beta}_r \v s_i^{(t)})^\T \,\Big|\, \v w_{i,c(i)}^{(t)}, \widehat{\v\mu}_{r}, \widehat{\m\Sigma}_{r}, \v\alpha_r', \m\Beta_r', \m\Kappa_r' \Bigr] - \widehat{\v\alpha}_r \widehat{\v\alpha}_r^\T.
	\label{eq:CEM_kappa}
\end{align}
The iteration is repeated with $(\v\alpha_r', \m\Beta_r', \m\Kappa_r') = (\widehat{\v\alpha}_r, \widehat{\m\Beta}_r, \widehat{\m\Kappa}_r)$ until convergence is effectively obtained, whereby a final estimate $(\widehat{\v\alpha}_r, \widehat{\m\Beta}_r, \widehat{\m\Kappa}_r)$ is returned.

\section{Conditional expectation of noisy observations}
\label{app:denoise}

Let the situation be as in Section\nobreakspace \ref {s:denoising} and consider the GLG model. 
The joint density of $(s,v)$ is found just as in the noise-free case in \textup {(\ref {eq:jointpdf_root})},
\begin{equation*}
	p(s, v | \mu, \sigma^2) 
	= p(v | s) p(s | \mu, \sigma^2) 
	= \frac{\exp\Bigl(-\frac{1}{2}\Bigl[\frac{v^2}{\exp(s) + \sigma_\varepsilon^2} + \frac{(s-\mu)^2}{\sigma^2}\Bigr]\Bigr)}{2\pi\sigma \sqrt{\exp(s) + \sigma_\varepsilon^2}}
\end{equation*}
and the marginal density of the wavelet with noise is
\begin{equation*}
	 p(v | \mu, \sigma^2) =
         \int_{-\infty}^\infty p(s, v | \mu,
         \sigma^2)\, \mathrm ds.
\end{equation*}
We do not have a closed form expression for this integral, but due to the form of the integrant we approximate the integral with the Gauss-Hermite quadrature rule, see e.g.\ \cite{Press:Teukolsky:Vetterling:Flannery:2002}.
The conditional density of $s$ given $v$ is
\begin{equation*}
	p(s | v, \mu, \sigma^2) = \frac{ p(s, v | \mu, \sigma^2) }{ p(v | \mu, \sigma^2) }
	= \frac{\exp\Bigl(-\frac{1}{2}\Bigl[\frac{v^2}{\exp(s) + \sigma_\varepsilon^2} + \frac{(s-\mu)^2}{\sigma^2}\Bigr]\Bigr)}{c(v | \mu, \sigma^2)\sqrt{\exp(s) + \sigma_\varepsilon^2}}
\end{equation*}
where $c(v | \mu, \sigma^2) = 2\pi\sigma p(v | \mu, \sigma^2)$. 
Furthermore, from well-known results about the bivariate normal distribution we obtain
\begin{equation*}
	\E[w | s, v, \v\theta] 
	=\frac{\cov[w, v | s, \v\theta]}{\var[v | s, \v\theta] } v 
	= \frac{ \var[w | s] }{ \var[v | s] } v
	= \frac{ \exp(s) }{ \exp(s) + \sigma_\varepsilon^2} v.
\end{equation*}
Hence
\begin{equation}\label{AHA}
	\E[w | v, \v\theta] 
	= \E\bigl[ \E[w | s, v, \v\theta] \big| v, \v\theta \bigr] 
	= v\, \E\biggl[ \frac{ \exp(s) }{ \exp(s) + \sigma_\varepsilon^2} \bigg| v, \v\theta \biggr]
\end{equation}
whereby we obtain \textup {(\ref {eq:freq_post_mean})}.

\bibliographystyle{plain}
\bibliography{literature}

\end{document}